\documentclass[12pt]{amsart}

\newcommand{\Z}{\mathcal{Z}}            
\newcommand{\W}{\mathcal{W}}            
\newcommand{\ID}{\mathbb{D}}            
\newcommand{\IN}{\mathbb{N}}            
\newcommand{\IC}{\mathbb{C}}            
\newcommand{\dbar}{\bar\partial}
\newcommand{\grad}{\nabla}
\newcommand{\lap}{\partial\bar\partial} 
\newcommand{\invlap}{\tilde \Delta}
\newcommand{\defeq}{\stackrel{\mathrm{def}}{=}}
\newcommand{\term}{\emph}
\newcommand{\st}{:}
\let\eps\epsilon

\newtheorem{theorem}{Theorem}[section]
\newtheorem{lemma}[theorem]{Lemma}
\newtheorem{proposition}[theorem]{Proposition}
\newtheorem{corollary}[theorem]{Corollary}
\newtheorem{Athm}{Theorem}

\theoremstyle{remark}
\newtheorem{remark}[theorem]{Remark}
\numberwithin{equation}{section}

\begin{document}

\title [Interpolating Sequences for the Bergman Space]
       {Interpolating Sequences for the Bergman Space\\
         and the $\bar\partial$-equation in Weighted $L^p$}
\author{Daniel H. Luecking}
\address{Department of Mathematical Sciences\\
         University of Arkansas\\
         Fayetteville, Arkansas 72701}
\email{luecking@uark.edu}
\subjclass{Primary 46E20}
\keywords{Bergman space, interpolating sequence, d-bar equation}

\begin{abstract}
The author showed that a sequence $\mathcal{Z}$ in the
unit disk is a zero sequence for the Bergman space $A^p$ if and only if
the weighted space $L^p(e^{pk_{\mathcal{Z}}}\,dA)$ contains a non-zero
(equivalently, zero-free) analytic function, where
\begin{equation*}
    k_\Z(z) = \sum_{a\in\Z} \frac{(1-|a|^2)^2}{|1 - \bar
    az|^2}\frac{|z|^2}{2}.
\end{equation*}
Here we show that $\mathcal{Z}$ is an interpolating sequence for $A^p$
if and only if it is separated in the hyperbolic metric and the
$\bar\partial$-equation
\begin{equation*}
    (1 - |z|^2)\bar\partial u = f
\end{equation*}
has a solution $u$ satisfying $\| u \|_{p,\Z} \le C\| f \|_{p,\Z}$ for
every $f \in L^p(e^{pk_{\mathcal{Z}}}\,dA)$, where $\| f \|_{p,\Z}$
denotes the $L^p(e^{pk_\Z})$ norm. This holds for all finite $p
\ge 1$, but also for $p < 1$ if $L^p(e^{pk_{\mathcal{Z}}}\,dA)$ is
suitably modified.

In addition, we show how this relates to a recent criterion for
weighted $\bar\partial$ estimates by J.~Ortega-Cerd\`a, and to
K.~Seip's criterion for interpolation.
\end{abstract}

\maketitle

\section{Introduction}\label{sec:intro}

Let $\ID$ denote the unit disk $\{ z \st |z| < 1 \}$ in the complex plane
$\IC$. For $0 < p < \infty$, let $L^p$ denote the usual Lebesgue space
consisting of those measurable functions $f$ on $\ID$ such that the norm
$\| f \|_p$ defined by $\| f \|_p^p = \int |f|^p \,dA$ is finite.
Here $dA$ denotes area measure. Let $A^p$ denote the subspace of $L^p$
consisting of analytic functions.

If $f$ is analytic in $\ID$ and not identically zero, let $\Z(f)$ denote
its zero sequence. This means that $\Z(f)$ consists of a listing (in
some order) of all the points where $f$ has a zero, each such point
being repeated as often as the multiplicity of the zero. We will view
$\Z(f)$ as a \term{set with multiplicity}: an equivalence class of such
sequences, where two sequences are equivalent if one is a reordering of
the other. We write $\Z(f)\subset \Z(g)$ if some representative of
$\Z(f)$ is a subsequence of a representative of $\Z(g)$.

A sequence $\Z = (a_n, n\in \IN)$ in $\ID$ is an \term{$A^p$ zero
sequence} if there is a non-zero function $f \in A^p$ such that $\Z(f) =
\Z$, where again we view $\Z$ as a set with multiplicity. Necessarily,
$\Z$ cannot have limit points in $\ID$ if it is a zero sequence.

We associate to a sequence $\Z$ in $\ID$ the normed sequence space
$l^p(\Z)$ consisting of all sequences $ c = (c_a, a\in \Z)$ such that
the norm $\| c \|_{p,\Z}$ defined by $\| c \|_{p,\Z}^p = \sum_{a\in\Z}
|c_a|^p(1 - |a|^2)^2$ is finite (this is a quasi-normed if $p < 1$). The
sequence $\Z$ is called an \term{interpolating sequence for $A^p$} if
the operator $f \mapsto (f(a) \st a \in\Z)$ takes $A^p$ continuously
\emph{onto} $l^p(\Z)$. An interpolating sequence is necessarily a zero
sequence for $A^p$, without repetition. A~criterion for $\Z$ to be
interpolating for $A^p$ was provided by K.~Seip in \cite{Sei94}, and
extended by A.~Schuster in \cite{Sch97}.

If $\Z$ is a sequence satisfying $\sum_{a\in \Z} (1 - |a|^2)^2 <
\infty$, define the function $k_{\Z}$ by
\begin{equation*}
    k_{\Z}(z) = \frac{|z|^2}{2}\sum_{a\in\Z} \frac{(1 - |a_n|^2)^2}{|1-
    \bar a_n z|^2}.
\end{equation*}
In \cite{Lue96} it was shown that $\Z$ is an $A^p$ zero sequence if and
only if there exists a zero-free function $f$, analytic in the unit
disk, such that $fe^{k_{\Z}}$ belongs to $L^p$. And this is equivalent
to there existing a non-zero function with this property. Let $L^p_{\Z}$
denote the space of all measurable functions $f$ such that $fe^{k_{\Z}}$
lies in $L^p$, let $\| f \|_{p,\Z}$ be the $L^p$ norm of this product,
and let $A^p_{\Z}$ denote the analytic functions in $L^p_{\Z}$. The main
result of \cite{Lue96}, then, is that $\Z$ is an $A^p$ zero sequence if
an only if $L^p_{\Z}$ contains a non-trivial analytic function, and so
$A^p_{\Z}$ is non-trivial.

The precise result shown in \cite{Lue96} was that if $\Psi_{\Z}$ denotes
the product
\begin{equation*}
    \Psi_{\Z}(z) = \prod \bar a_n \frac{a_n - z}{1 - \bar a_n z}
    \exp\left( 1 - \bar a_n \frac{a_n - z}{1 - \bar a_n z} \right),
\end{equation*}
then $\Z$ is an $A^p$ zero set if and only if division by $\Psi_{\Z}$ is
an isomorphism from the subspace $I^p_{\Z} = \{ f\in A^p \st
\Z(f) \supset \Z\}$ onto $A^p_{\Z}$. (This is
strictly correct only if $0 \notin \Z$, but a correct version is easily
supplied in that case.)

The starting point of our investigations was to try to relate properties
of the sequence $\Z$ to properties of $k_{\Z}$ and in particular to
properties of $A^p_{\Z}$ and $L^p_{\Z}$. Since $\Z$ being a zero
sequence for $A^p$ is equivalent to the property that $L^p_{\Z}$
contains non-trivial analytic functions, one might ask what further
properties of $L^p_{\Z}$ are equivalent to $\Z$ being an interpolating
sequence. It is natural to conjecture that $\Z$ is an interpolating
sequence if and only if there are certain $L^p_{\Z}$ bounds on the
solution of the $\dbar$-equation. Part of the reason this is natural is
because one direction is rather routine to prove ($\dbar$-estimates
imply interpolation).

Let $w$ denote a weight function on $\ID$, that is, $w$ is measurable,
not necessarily integrable, and $w(z) > 0$ for all
$z\in \ID$. Let $L^p(w)$ denote the weighted space of measurable
functions with norm $\| f \|_{p,w}^p = \int |f|^p w\,dA < \infty$. The
$\dbar$ and $\partial$ operators are defined by
\begin{align*}
  \dbar u &= \frac{\partial u}{\partial \bar z} = \frac{\partial
  u}{\partial x} - \frac{1}{i} \frac{\partial u}{\partial y} \\
  \partial u &= \frac{\partial u}{\partial z} = \frac{\partial
  u}{\partial x} + \frac{1}{i} \frac{\partial u}{\partial y}.
\end{align*}
We will consider the equation
\begin{equation}\label{eq:dbar}
    (1 - |z|^2)\dbar u = f
\end{equation}
and try to determine, for certain weights $w$, when there exists a
constant $C$ such that it can always be solved with $\| u \|_{p,w} \le C
\| f \|_{p,w}$. The above mentioned conjecture deals with $L^p_{\Z}$,
which is just the case $w = e^{pk_{\Z}}$. The main result of this
article is the following validation of that conjecture.

\begin{Athm}\label{thm:mainA}
  A sequence $\Z$ in the unit disk is interpolating for $A^p$, $p \ge
  1$, if and only if it is separated in the hyperbolic metric and there
  is a constant $C$ such that for every $f \in L^p_{\Z}$, there is a
  solution of equation~\eqref{eq:dbar} satisfying $\| u
  \|_{p,\Z} \le C\| f \|_{p,\Z}$.
\end{Athm}

In section~\ref{sec:background} we present the background on $k_{\Z}$
and other topics, and include some of the basic lemmas needed for
proving the main results. In section~\ref{sec:main1}, we prove one
direction of the main result, that the existence of solutions of the
$\dbar$-equation with the given bounds implies interpolation. Then, in
section~\ref{sec:main2}, we prove the converse. In section~\ref{sec:plt1}
we show how the main result may be modified to provide a version that
holds for $p < 1$. All this relies only on simple properties of
interpolating sequences and rather easy properties of the connection
between zero sequences and $k_{\Z}$ set out in \cite{Lue96} and
\cite{Lue99c}. In fact, this author believes the result itself is not as
important as the variety of techniques used in proving it.

In the final two sections we explore the connection between K. Seip's
criterion for interpolation (a certain uniform upper bound on the
density of the sequence $\Z$) and a recent criterion by
J.~Ortega-Cerd\`a for the solution of \eqref{eq:dbar} in weighted $L^p$
with certain types of weights. In section~\ref{sec:seip} we show that
Seip's criterion translates to an upper bound on certain averages of
$\lap k_{\Z}$ over disks of a fixed pseudo-hyperbolic radius. In
section~\ref{sec:ortega}, we observe that Ortega-Cerd\`a's criterion
(for solutions of the $\dbar$ equation) and this form of Seip's
criterion (for interpolation) are remarkably similar. While
Ortega-Cerd\`a's criterion is stated in a context that, on the face of
it, excludes $e^{k_{\Z}}$ as the weight, we can show that the weight is
equivalent one that does satisfy Ortega-Cerda\'a's criterion. This
produces an alternate proof of the main result, showing that Seip's
criterion also implies solutions of the $\dbar$ equation with bounds in
$L^p_\Z$.

Many of the auxiliary results we will prove are already known for $A^p$,
and in some cases much stronger results are known. However, we will use
the weaker versions and provide proofs because, (i)~the proofs are new;
(ii)~they rely only on quite simple properties of zero/interpolating
sequences; and (iii)~they extend to weighted Bergman spaces where the
stronger results do not. In fact, the above Theorem remains valid
(appropriately modified) for all the usual weighted spaces
$A^{p,\alpha}$, $\alpha> -1$. We will always try to make it clear when
we are reproving a known result and to state what settings our proof
extends to. In most cases the proof will extend almost without change to
the weighted spaces.

I would like to gratefully acknowledge the help of J.~Ortega-Cerd\`a,
who showed me the methods of section~\ref{sec:ortega} after reading a
preliminary version of this paper.

\section{Background and auxilliary results} \label{sec:background}

Let $\phi(z)$ denote a real-valued function on $\ID$ and
denote by $L^p_\phi$ the space $L^p(e^{p\phi})$, that is
\begin{equation*}
    L^p_\phi = \left\{ f \st \int_{\ID} \left| f(z) e^{\phi}\right|^p
    \,dA(z) < \infty \right\}.
\end{equation*}
The norm $\| f \|_{p,\phi}$ will denote the $p$-th root of the integral
in the above definition. We will always assume that $e^{p\phi}$ is
locally integrable, so that all continuous functions with compact
support in $\ID$ belong to $L^p_\phi$. The space $L^p_{\Z}$ is just
$L^p_\phi$ with $\phi(z) = k_{\Z}(z)$.

We will say that we \term{can solve the $\dbar$-equation with bounds} in
$L^p_{\phi}$ if there is a constant $C > 0$ such that equation
\eqref{eq:dbar} has, for all $f\in L^p_\phi$, a solution $u \in
L^p_\phi$ satisfying $\| u \|_{p,\phi} \le C\| f \|_{p,\phi}$. We
understand the equation to hold either in the sense of distributions, or
in the sense that it holds classically for all $f$ which are continuous
with compact support. The hypotheses on $\phi$ make such functions dense
in $L^p_{\phi}$.

Recall that the pseudo-hyperbolic metric $\psi$ is defined by
\begin{equation*}
    \psi(z,w) = \left| \frac{z - w}{1 - \bar wz} \right|
\end{equation*}
For $a \in \ID$ we denote by $M_a$ the M\"obius transformation
\begin{equation*}
    M_a(z) = \frac{a-z}{1 -\bar az}.
\end{equation*}
Of course, the metric $\psi$ is invariant under M\"obius transformations,
and is locally equivalent to the Euclidean metric in $\ID$.

Define $D(z,r)$ to be the ball in the metric $\psi$ with center $z$ and
radius $r$. Note that $D(z,r)$ is the Euclidean disk centered at $z(1-
r^2)/(1 - r^2|z|^2)$ with radius $r(1 - |z|^2)/(1 - r^2|z|^2)$.

We will have need of the following lemma which provides a test, often
called the Schur criterion, for the boundedness in $L^p$ of integral
operators.

\begin{lemma}\label{lem:schur} \textup{(Schur criterion)}
  Let $K(z,w)$ be a measurable kernel on measure spaces $(Z,\nu) \times
  (W,\mu)$ and let $p > 1$. If there exist measurable functions $h_1(w)
  > 0$ and $h_2(z) > 0$, and constants $C_1$ nd $C_2$ such that
  \begin{align}
    \int |K(z,w)| h_1(w)^{p'} \,d\mu(w) &\le C_1 h_2(z)^{p'}\qquad w\in W
      \label{eq:schur1}\\
\intertext{and}
    \int |K(z,w)| h_2(z)^p \,d\nu(z) &\le C_2 h_1(z)^p \qquad z\in Z\label{eq:schur2}
  \end{align}
  where $p'=p/(p - 1)$, then the operator $T$ defined by $Tf(z) = \int
  K(z,w) f(w) \,d\mu(w)$ is bounded from $L^p(W,\mu)$ to $L^p(Z,\nu)$
  and $\|T\| \le C_1^{1/p'}C_2^{1/p}$. Moreover, $T$ will be bounded
  from $L^1(W,\mu)$ to $L^1(Z,\nu)$, with $\| T \| \le C$ if
  \begin{equation*}
    \sup_{w\in W} \int |K(z,w)| \,d\nu(z) \le C.
  \end{equation*}
\end{lemma}

What we will mostly need is the following consequence.

\begin{lemma}\label{lem:boundedkernel}
 If $p \ge 1$, the kernels
  \begin{gather*}
    K(z,w) =  \frac{(1 - |z|^2)^a(1 - |w|^2)^b}{|1 - \bar w
    z|^{a+b+2}}\\
    B(z,w) =  \frac{(1 - |z|^2)^a(1 - |w|^2)^b}{|z-w| \, |1 - \bar w
    z|^{a+b+1}}
  \end{gather*}
  define bounded operators on $L^p(dA)$ \textup(integration with respect to
  $w$\textup) provided $a > -1/p$ and $b > -1/p'$. In case $p = 1$ the
  condition on $b$ means $ b > 0$.
\end{lemma}

\begin{proof}
It is well known (see for example \cite{Zhu90}) that if $-1 < \beta < M$
then
\begin{equation}\label{eq:Forelli-Rudin}
  \int \frac{(1- |w|^2)^\beta}{|1 - \bar wz|^{M+2}} \,dA(w) \le C(1 -
  |z|^2)^{\beta - M}.
\end{equation}
We observe that there is a similar estimate
\begin{equation}\label{eq:FR2}
  \int \frac{(1- |w|^2)^\beta}{|z - w|\,|1 - \bar wz|^{M+1}} \,dA(w) \le C(1 -
  |z|^2)^{\beta - M}.
\end{equation}
This can be seen by writing the integral as the sum of two integrals,
one over the set where $|\psi(z,w)| \ge \eta$ (for some $0 < \eta < 1$)
and one over the complement. The first integral can be estimated as at
most $1/\eta$ times the integral in \eqref{eq:Forelli-Rudin}. In the
second integral we can use the fact that, for $w \in D(z,\eta)$, we have
\begin{equation*}
 (1 - |z|^2)/C_\eta \le 1 - |w|^2 \le C_\eta (1 - |z|^2) \quad
 \hbox{and}\quad |1 - \bar wz| > (1 - |z|^2)/C_\eta
\end{equation*}
for some constant $C_\eta$ depending only on $\eta$. This gives
\begin{equation*}
  \int_{D(z,\eta)} \frac{(1- |w|^2)^\beta}{|z - w|\,|1 - \bar wz|^{M+1}}
  \,dA(w) \le [C_\eta(1 - |z|^2)]^{\beta - M - 1} \int_{D(z,\eta)}
  \frac{1}{|z-w|} \,dA(w).
\end{equation*}
Finally, the integral on the right side above is equivalent to the
Euclidean radius of $D(z,\eta)$, which is bounded by $C_\eta (1 -
|z|^2)$.

Let $p > 1$. In Lemma~\ref{lem:schur}, let $h_1(z) = (1 - |z|^2
)^{-\alpha}$ and let $h_2 = h_1$. Letting $M = a+b$, apply
inequality~\eqref{eq:Forelli-Rudin} or \eqref{eq:FR2} twice, once with
$\beta = b-p'\alpha$ and once with $\beta=a-p\alpha$, then both
\eqref{eq:schur1} and \eqref{eq:schur2} will be satisfied if we can find
$\alpha$ satisfying the following simultaneous inequalities.
\begin{gather*}
  -1 < b - \alpha p' < a + b\\
  -1 < a - \alpha p < a + b
\end{gather*}
These are equivalent to
\begin{gather*}
   -a/p' < \alpha < (1 + b)/p'\\
   -b/p < \alpha < (1 + a)/p
\end{gather*}
The existence of $\alpha$ satisfying all four inequalities is equivalent
to
$$
  \max\{-a/p',-b/p\} < \min\{ (1 + b)/p', (1 + a)/p \}.
$$
This gives four inequalites, two of which are equivalent to each other
and implied by the other two, reducing the requirement to
\begin{align*}
  -a/p' &< 1/p  + a/p\\
  -b/p  &< 1/p' + b/p'
\end{align*}
which are equivalent to the two requirements in the statement of this
lemma.

The $p=1$ case is immediate from \eqref{eq:Forelli-Rudin} or
\eqref{eq:FR2} with $\beta=b$ and $M=a+b$.
\qed\end{proof}

We need some estimates for the maximum size of an analytic function at
the origin, given a bound on its norm plus the condition that it vanish
exactly on an interpolating sequence $\Z$.

A sequence is said to be \term{separated} in the pseudo-hyperbolic metric if
there is a positive lower bound on $\psi(z,w)$ for $z \ne w$ in $\Z$.
The phrase \term{uniformly discrete} is also used with the same meaning.
The largest lower bound will be referred to as the \term{separation
constant}.

Now if $\Z$ is an interpolating sequence then it is a zero sequence for
$A^p$ and it is also separated in the pseudo-hyperbolic metric. Either
condition forces $\sum_{a\in \Z} (1 - |a|^2)^2 < \infty$. Define for $a
\in \ID$,
\begin{equation}
  E_a(z) =
    \begin{cases}
      \dfrac{\bar a}{|a|} \dfrac{a - z}{1 - \bar az}
        \exp\left[1 - \bar a\dfrac {a - z} {1 - \bar az} - (1 -
        |a|^2)/2\right] & a \ne 0\\[10pt]
      ze^{1/2} & a = 0
    \end{cases}
\end{equation}
and then let $\Psi_{\Z}(z) = \prod_{a\in\Z} E_a(z)$.
This is almost the same function as in \cite{Lue96}, but it has been
normalized by dividing each factor by $|a|\exp[(1 - |a|^2)/2]$, and an
oversight corrected so that $0 \in \Z$ is properly accounted for.
This function $\Psi_{\Z}$ is analytic in $\ID$ and vanishes only on
$\Z$. The odd form of the factor for zeros at the origin is for
consistency with the other factors and also so that \eqref{eq:Psi.kZ}
below holds.

Now let
\begin{equation}
  \sigma_{\Z}(z) = \prod_{a\in Z} \left| \frac{a - z}{1 - \bar az}
  \right|\exp\left[ \frac{1}{2}\left( 1 - \left| \frac{a - z}{1 - \bar az}
  \right|^2 \right) \right].
\end{equation}
Here no special provision needs to be made for $a = 0 \in \Z$.
What was shown in \cite{Lue96} is that there is a constant $C_p$
depending only on $p$ such that if $f \in A^p$ vanishes on $\Z$ then
\begin{equation}
    \| f/\sigma_{\Z} \|_p \le C_p \| f \|_p.
\end{equation}
We remark here the rather surprising fact that $C_p$ \emph{does not
depend on $\Z$}. This follows from the proof in \cite{Lue96}, although
that paper does not mention it. Since
\begin{equation}\label{eq:Psi.kZ}
    \left| \Psi_{\Z} \right| = \sigma_{\Z}e^{k_{\Z}},
\end{equation}
it follows that division by $\Psi_{\Z}$ takes the subspace $I^p_{\Z}$ of
functions in $A^p$ that vanish on $\Z$ into $A^p_{\Z}$ boundedly, with a
bound independent of $\Z$. Moreover, since $\sigma_{\Z} \le 1$, it also
follows that multiplication by $\Psi_{\Z}$ is a bounded from $A^p_{\Z}$
onto $I^p_{\Z}$, with norm at most $1$. These results are actually the
easiest part of \cite{Lue96}, requiring only an integration of Jensen's
formula followed by the inequality of geometric and arithmetic means,
followed by the Schur criterion (Lemma~\ref{lem:schur} above).

Now observe the simple fact that if $\W\subset\Z$ (taking multiplicities
into account) then $A^p_{\Z} \subset A^p_{\W}$ with the inclusion map
bounded by 1 in norm. This is because the weights $\exp k_{\Z}$ increase
with $\Z$. The composition of these three operators:  division by
$\Psi_{\Z}$, this inclusion, and multiplication by $\Psi_{\W}$ is
therefore a bounded map from $I^p_{\Z}$ into $I^p_{\W}$, with a bound
that depends only on $p$. This operator is obviously the same as
division by $\Psi_{\Z\setminus \W}$. This makes the following lemma
almost obvious.

\begin{lemma}\label{lem:removezeros}
  If $f\in A^p$ satisfies $\| f \|_p = 1$, $|f(0)| \ge \delta$ and $\W
  \subset \Z(f)$, then there is a function $g\in A^p$ such that
  $\W=\Z(g)$, satisfying $\| g \|_p = 1$ and $|g(0)| > \delta/C_p$,
  where $C_p$ is a positive constant depending only on $p$.
\end{lemma}

We note that this lemma actually holds with $C_p = 1$ on $A^p$ by
\cite{DKS93a}. While our result is weaker, it suffices for our needs.
Moreover, our result extends without change to the weighted case
$A^{p,\alpha}$ with a constant $C_{p,\alpha}$ depending only on $p$ and
$\alpha$ and in fact to all the spaces covered in \cite{Lue96}, with the
constant depending only on the space.

\begin{proof}
Apply the preceding discussion with $\Z$ being equal to the zero set of
$f$, producing $h = f/\Psi_{\Z'}$ where $\Z' = \Z \setminus \W$. Then
$\| h \|_p \le C_p$ with a constant $C_p$ depending only on $p$. Now
\begin{equation*}
    \log |\Psi_{\Z'}(0)| = \sum_{a\in \Z'} \log |a| + \frac{1}{2}(1 -
    |a|^2)
\end{equation*}
is easily seen to be negative, so $|\Psi_{\Z'}(0)| < 1$. Therefore
$|h(0)| > |f(0)| > \delta$ and so $g = h/\| h \|_p$ has the required
properties.
\qed\end{proof}

The \term{interpolation constant} for an interpolating sequence $\Z$ is
the least number $M$ such that if $(c_a \st a\in\Z)$ is a sequence in
the unit ball of $l^p_{\Z}$ then there exists a function $f\in A^p$ such
that $f(a) = c_a$ for all $a\in \Z$ and $\| f \|_p \le M$. The existence
of $M$ is a consequence of the open mapping theorem since the operation
of evaluation on $\Z$ is bounded from $A^p$ \emph{onto} $l^p_{\Z}$.

The following seems to be well-known, but we include a proof for
completeness, and because it is more elementary than some proofs that
the author has seen, and because a published proof was hard to find.

\begin{lemma}\label{lem:addpoint}
  Let $\Z$ be an interpolation sequence for $A^p$ with interpolation
  constant $M$. Let $0 < \eta < 1$, and let $a_0$ be a point of $\ID$
  satisfying $\psi(a_0,a) > \eta$ for all $a\in\Z$. Then $\Z \cup \{ a_0
  \}$ is also an interolating sequence for $A^p$ with interpolation
  constant less than $C/\eta$, where $C$ depends only on $M$ and $p$.
\end{lemma}

\begin{proof}
Let $\Z$ be interpolating and let $\Z' = \Z \cup \{ a_0 \}$. Let $c \in
l^p_{\Z'}$ with $\| c \|_{p,\Z'} = 1$. We seek a function $f$ in $A^p$
with $f(a) = c_a$ for all $a \in \Z'$, together with an estimate of its
norm. Since a M\"obius transform of an interpolating sequence is also
interpolating with the same interpolation constant we may assume $a_0 =
0$ and therefore $|a| > \eta$ for all $a \in \Z$. First obtain a
function $g \in A^p$ with $g(a) = (c_a - c_0)/a$ for $a \in Z$. The norm
of $\left( (c_a - c_0) \st a \in \Z' \right)$ is easily estimated in terms
$p$ and $\sum (1 - |a|^2)^2$, and the latter depends only on the
separation constant, which can be seen to depend only on $M$ and $p$.
Therefore, the norm of $\left( (c_a-c_0)/a \right)$ is at most $C'/\eta$,
and so such a $g$ can be found with $\| g \|_p \le MC'/\eta$. Then
let $f(z) = zg(z) + c_0$ to get the required interpolation with norm at
most $MC'/\eta + 1 \le C/\eta$.
\qed\end{proof}

Combining this lemma with the previous one we obtain the result we will
later need.

\begin{corollary}\label{cor:lowerestimate}
  Let $\Z$ is an interpolating sequence for $A^p$ with interpolation
  constant $M$, let $0 < \eta < 1/2$ and let $\W = \Z \setminus
  D(0,\eta)$. There is a constant $\delta > 0$ depending only on $p$,
  $M$ and $\eta$ such that there exists a function $f \in A^p$ with $\|
  f \|_p \le 1$, $\Z(f) = \W$ and $|f(0)| > \delta$.
\end{corollary}

\begin{proof}
The sequence $\W$ will be interpolating with constant at most $M$, the
sequence $\W' = \{ 0 \} \cup \W$ will be interpolating with constant at
most $C/\eta$, where $C$ depends only on $M$ and $p$. Select a function
$g\in A^p$ with at most this norm that is $1$ at the origin and $0$ on
the rest of $\W'$. Let $h = g/\| g \|_p$ so that $|h(0)| \ge \eta/C$.
Replace $h$, if necessary, by a function $f$ that vanishes only on
$\W$. This can be done, by Lemma~\ref{lem:removezeros},
with a function $f$ satisfying $\| f \|_p = 1$ and $|f(0)| \ge
|h(0)|/C_p$. This proves the corollary with $\delta \ge \eta/(C_pC)$.
\qed\end{proof}

The following stability result for interpolation sequences was proved in
\cite{Lue85a}:

\begin{proposition}\label{prop:interpolationperturb}
  If $\Z = \{ a_n \}$ is an interpolating sequence for $A^p$ with
  interpolation constant $M$ then there exists a constant $\delta > 0$
  depending only on $M$ and $p$ such that any sequence $\W = \{ b_n \}$
  satisfying $\psi(a_n, b_n) < \delta$ is also interpolating for $A^p$
  with interpolation constant less than $2M$.
\end{proposition}

Actually, the constant for the new sequence $\W$ can can be taken to be
arbitrarily close to $M$ if $\delta$ is made sufficiently small, but we
will not need that fact.

Another necessary result, the following was mostly proved in
\cite{Lue99c}. If $0 < \lambda < 1$ and $a \in \ID$ is nonzero, let
$p_\lambda(a)$ denote the point with the same argument as $a$ satisfying
$1 - |p_\lambda(a)|^2 = \lambda(1 - |a|^2)$. Let $p_\lambda(0) = 0$.

\begin{proposition}\label{prop:zeroperturb}
  Suppose $0 \notin Z$ and the sequence $\Z \cup \{ 0 \}$ is separated
  in the pseudo-hyperbolic metric with separation constant $\eta$.  Let
  $0 < \lambda < 1$ and let $p_\lambda(\Z)$ denote the set of all
  $p_\lambda(a)$ for $a\in\Z$. If $\Z$ is a zero sequence for $A^p$ then
  $p_\lambda(\Z)$ is a zero set for $A^{p/\lambda}$. Moreover, if $f \in
  A^p$ with $\Z= \Z(f)$, then there exists an $h \in A^{p/\lambda}$ and
  a constant $C > 0$ depending only on $\lambda$ and $\eta$ such that
  $\| h \|_{p/\lambda} \le \| f \|_p$ and $|h(0)| > |f(0)|^\lambda/C$.
\end{proposition}

The last part of this above proposition was not mentioned in
\cite{Lue99c}, but follows from its method of proof. For completeness,
we outline the proof here.

\begin{proof}
For simplicity of notation, let $\Z'=p_\lambda(\Z)$ and for each
$a\in\Z$ let $a'=p_\lambda(a)$ Given an $f\in A^p$ with zero sequence
$\Z$, it follows from \cite{Lue96} that $g = f/\Psi_\Z$ is a nowhere
zero function in $A^p_\Z$. In \cite{Lue99c}, the hypotheses were shown
to imply that the following sum converges to a harmonic function
\begin{equation*}
  u(z) =
  \sum_{a\in\Z} (1 - |a|^2) \left[ \frac {1 - |a'|^2 |z|^2}
  { |1 - \bar a' z |^2 } - \frac{1 - |a|^2 |z|^2}
  { |1 - \bar az|^2} \right]
\end{equation*}
It was shown, moreover, that
\begin{equation*}
   k_{\Z'}(z) - \lambda k_\Z(z) -  u(z)
\end{equation*}
is bounded, with the bound $c$ depending only on $\lambda$ and the
separation constant of $\Z$. Let $v$ be harmonic in $\ID$ with $u+ iv$
analytic and put $h_1 =  g^\lambda e^{-u-iv}$, then
\begin{equation*}
  \left| h_1 e^{k_{\Z'}}\right|^{p/\lambda} \le e^c \left| g
  e^{k_\Z} \right|^p
\end{equation*}
Since $u(0) = 0$, we see that $|h_1(0)| = |g(0)|^\lambda$. Now let $h=
h_1\Psi_{\Z'}$. Then we have
\begin{equation*}
  \| h \|_{p/\lambda} \le \| h_1 \|_{p/\lambda,\Z'} \le C \| g
  \|_{p,\Z} \le C'\| f \|,
\end{equation*}
and
\begin{equation*}
  |h(0)| = |f(0)|^\lambda \frac{ | \Psi_{\Z'} (0) |}
  { | \Psi_{\Z} (0) |^\lambda }.
\end{equation*}
Normalizing $h$ to have the same norm as $f$ gives the required result,
since $|\Psi_\Z(0)|$ and $|\Psi_{\Z'}(0)|$ can be estimated in terms of
the separation constant $\eta$.
\qed\end{proof}

\begin{proposition}\label{prop:growth}
  Let $\Z$ be any zero set for $A^p$, $p > 0$, and suppose $0
  \notin \Z$. Then there exists a solution $f^*$ of the extremal
  problem: maximize $|f(0)|$ for $f \in A^p$ subject to $\| f \|_p
  \le 1$ and $\Z(f) = \Z$. This solution satisfies
  \begin{equation}\label{eq:harm.eval}
    \int |f^*|^p u \,dA = u(0)
  \end{equation}
  for all bounded harmonic $u$ on $\ID$.
  As a consequence of \eqref{eq:harm.eval}, there is a
  constant $C$ such that
  \begin{equation}\label{eq:growth}
    |f^*(z)|^p (1 - |z|^2) \le C \quad\mbox{for all $z \in \ID$,}
  \end{equation}
  and $g = f^*/\Psi_\Z$ satisfies a similar growth condition
  \begin{equation}\label{eq:growth2}
    |g(z)|^p e^{pk_{\Z}(z)} (1 - |z|^2) \le C'
        \quad\mbox{for all $z \in \ID$.}
  \end{equation}
\end{proposition}

Note: In \cite{DKS93a} a similar extremal problem is considered, but the
constraint is weakened to require only $\Z \subset \Z(f)$. Then the
solution is \emph{proved} to satisfy $\Z(f^*) = \Z$. Here we maximize
only over the set of functions already having the given zero set. This
has the advantage of producing a result that remains valid in the
weighted spaces $A^{p,\alpha}$.

\begin{proof}
The solution exists by the following normal families argument. We
observe that the unit ball of $A^p$ is a normal family. Let $B$ be the
supremum of $|f(0)|$ over the unit ball of $A^p$. Let $f_n$ be a
sequence in the unit ball with $|f_n(0)| \to B$ and let $f^*$ be any
normal limit point of the sequence. Clearly $|f^*(0)| = B$ and $\|
f^* \|_p \le 1$. If the norm were not equal to $1$ we could increase
$|f^*(0)|$ by dividing $f^*$ by its norm, contradicting the definition
of $B$.

The dual problem (minimize $\| f \|_p$ subject to $|f(0)| = B$) has the
same solutions. Let $u$ be any bounded harmonic function on $\ID$ and $v$
its harmonic conjugate. Then for all real $t$, $f^* e^{t(u-u(0) + iv)}$
are among the candidates for this minimization, and so
\begin{equation*}
  \int |f^* e^{t(u-u(0) + iv)}|^p \,dA
\end{equation*}
has a minimum at $t = 0$. Setting the derivative at $t=0$ to zero gives
\begin{equation*}
  \int |f^*|^p  (u - u(0))\,dA = 0
\end{equation*}
which implies \eqref{eq:harm.eval}.

To get \eqref{eq:growth} we apply \eqref{eq:harm.eval} to a sequence of
bounded harmonic functions converging uniformly on compact sets to the
Poisson kernel $P_\zeta(z) = (1 - |z|^2)/|\zeta - z|^2$, $|\zeta| = 1$,
to get, after Fatou's lemma:
\begin{equation*}
  \int |f^*|^p  P_\zeta \,dA \le 1.
\end{equation*}
Let $z\in \ID$ with $z/|z| = \zeta$ and consider the restriction of the
above integral to the Euclidean disk $D_z=\{ w \st |w-z| < (1 - |z|)/2 \}$.
Then the above inequality combined with a routine estimate of $P_\zeta$
on $D_z$ gives:
\begin{equation*}
  \frac{1}{|D_z|} \int_{D_z} |f^*|^p  \le \frac{C}{1-|z|^2},
\end{equation*}
where $|D_z|$ denotes the area of $D_z$.

Finally, we have only to observe that, by subharmonicity, $|f^*(z)|^p$
is less than the above average. This gives inequality \eqref{eq:growth}.
Thus, $f^*$ belongs to one of the growth spaces $A(\infty, \infty, 1/p)$
considered in \cite{Lue96}. The results of that paper show that
inequality \eqref{eq:growth2} follows from \eqref{eq:growth}.
\qed\end{proof}

Note: The inequality \eqref{eq:growth2} (for some zero-free $g \in
A^p_\Z$) is the one we will need later. One could obtain it directly by
maximizing $|g(0)|$ among zero-free functions with norm 1 in $A_\Z^p$,
and making use of the fact that $\log |g| + k_\Z$ is subharmonic.

In the case of weighted spaces $A^{p,\alpha}$, we find that
$|f^*(z)|^p(1 - |z|^2)^{\alpha+1}$ is bounded. We also note that we can
interpolate: $f^*(z)(1 - |z|^2)^{1/p}$ belongs to both $L^\infty$ and
$L^p$ for the measure $(1 - |z|^2)^{-1}$ and therefore it belongs
simultaneously to $L^q((1 - |z|^2)^{-1}dm(z))$ for all $q \ge p$. That
is, if $\Z$ is an $A^p$ zero set, it is also an $A^{q,q/p - 1}$ zero
set for all $q > p$. A considerably stronger version of this observation
was proved by K.~Seip in \cite{Sei95}.

\section{Interpolating a given sequence}\label{sec:main1}

Here we prove the first part of the main theorem.

\begin{theorem}\label{thm:easydirection}
Assume that $\Z$ is separated in the pseudo-hyperbolic metric and that
we can solve the $\dbar$-equation with bounds in $L^p_{\Z}$, $p \ge 1$.
Then $\Z$ is an interpolating sequence for $A^p$.
\end{theorem}

\begin{proof}
Let $(c_a \st a\in \Z)$ be a sequence in $l^p_{\Z}$, and let $D_a =
D(a,2\eta)$, $a\in \Z$, be disjoint disks with equal pseudo-hyperbolic
radius. Let $D'_a = D(a,\eta)$.  Since the Euclidean radius of $D_a$ is
on the order of $C(1 - |a|)\eta$, we can find $C^1$ functions $\beta_a$
with support in $D_a$ satisfying
\begin{gather*}
   0 \le \beta_a \le 1\\
   \beta_a(z) = 1 \text{ for } z\in D'_a\\
   |\grad \beta_a(z)|(1 - |a|) \le C
\end{gather*}
It is clear that
\begin{equation*}
    g(z) = \sum_{a\in\Z} c_a \beta_a(z)
\end{equation*}
satisfies the interpolation $g(a) = c_a$, and moreover it belongs to
$L^p$ with norm $\| g \|_p \le C\| c \|_{p,\Z}$. We correct $g$ by
putting $f = g - u\Psi_Z$ and try to determine $u$ so that $f$ is
analytic. This requires $\dbar u = \dbar g/\Psi_{\Z}$ or
\begin{equation*}
    (1 - |z|^2)\dbar u = \frac{(1 - |z|^2)\dbar g}{\Psi_{\Z}}
\end{equation*}
We will show momentarily that the right hand side of this belongs to
$L^p_{\Z}$, with norm at most $C\| c \|_{p,\Z}$. Assuming this for the
moment, we then have a solution $u$ belonging to the same space with a
similar norm estimate. Because $|\Psi_{\Z}| \le e^{k_{\Z}}$, this means
$u\Psi_{\Z}$ and therefore $f = g + u\Psi_{\Z}$ belongs to $L^p$.
Moreover, it is analytic and satisfies $f(a) = c_a$, as required.

Since $|\Psi_{\Z}| = \sigma_{\Z}e^{k_{\Z}}$, we need only estimate
the quotient $(1 - |z|^2)\dbar g/\sigma_{\Z}$. Now $\dbar g$ is zero
inside each $D'_a$ as well as outside the union of the $D_a$. Moreover,
in each annular region $D_a \setminus D'_a$ it is bounded by $C c_a/(1 -
|a|^2)$. It therefore suffices to show that $\sigma_{\Z}(z) \ge \delta >
0$ for $z$ in each $D_a \setminus D'_a$. We do this by estimating from
above
\begin{equation*}
-2\log \sigma_{\Z}(z) =
    \sum_{a\in Z} \left[ \log \left| \frac{a - z}{1 - \bar az}
    \right|^{-2}
    - \left( 1 - \left| \frac{a - z}{1 - \bar az} \right|^2 \right) \right].
\end{equation*}
Now, if $1 > x > \eta > 0$,
\begin{equation*}
\log 1/x = \sum_{n = 1}^\infty \frac{(1-x)^n}{n} = 1 - x + (1 -
x)^2\sum_{k=0}^\infty \frac{(1 - \eta)^k}{k + 2} \le 1 - x + C_\eta (1 -
x)^2
\end{equation*}
so, with $x = |(z-a)/(1 - \bar az)|^2$, we have
\begin{equation*}
    -2\log \sigma_{\Z}(z) \le C_{\eta} \sum_{a \in \Z} \left( 1 - \left|
        \frac{a - z}{1 - \bar az} \right|^2 \right)^2
\end{equation*}
as long as $|a-z|/|1 - \bar az| > \eta$ for all $a\in Z$, which is true
if $z \in D_a \setminus D'_a$. This last sum is finite, with a bound
that can be estimated solely in terms of the separation constant of the
sequence
\begin{equation*}
    \left\{ \frac{a - z}{1 - \bar az} \st a \in \Z \right\}.
\end{equation*}
But that separation constant is the same as that of $\Z$, so the sum has
an upper bound independent of $z$.

\qed\end{proof}

\section{A solution operator for the $\dbar$-equation}\label{sec:main2}

We will construct a kernel for the solution of the $\dbar$-equation,
\eqref{eq:dbar}. Our method will be to construct local solutions and
patch them together with the aid of a family of analytic functions
$g_a$, $a\in \ID$, with special properties. The properties we need are:
$|g_a(z)e^{k_\Z(z)}|^p(1 - |M_a(z)|^2)^{1-\eps} \le C$, and
$|g_a(z)|e^{k_\Z}(z) > \delta > 0$ in $D(a, \eta)$ for some $C$, $\eps$,
$\delta$ and $\eta$ independent of $a$. We obtain this by using the
results of section \ref{sec:background}. These results imply that any
sequence $\Z$ that is interpolating for $A^p$ is a separated zero
sequence for $A^q$, some $q > p$, and moreover we have norm estimates on
functions that vanish on $\Z$ but are bounded below at $0$.

We will first do the construction of $g_a$ for $a = 0$ and then invoke
invariance of the hypotheses under M\"obius transformations of $\Z$ to
obtain such a function for each $a \in \ID$.

A rough outline of the construction of $g_0$ is the following. Take the
interpolating sequence $\Z$, assume it stays away from the origin, and
perturb each $a$ to $a'$ which has the same argument as $a$ but $1 -|a'|^2$
is \emph{increased} by a factor $1/\lambda$ with $\lambda < 1$. Make
$\lambda$ so close enough to $1$ that the resulting sequence is still
interpolating for $A^p$ (Proposition~\ref{prop:interpolationperturb}).
Then perturb it \emph{back} to where it was and it becomes a zero
sequence for $A^q$, $q = p/\lambda$, with the special properties of
Proposition~\ref{prop:zeroperturb}, namely the existence of a function
with norm 1 with zeros only on $\Z$ having a lower bound at the origin
depending on $p$ and the interpolation constant of $\Z$. Finally we
apply proposition~\ref{prop:growth}, possibly increasing the value at
the origin and forcing the appropriate growth estimates, namely
$|g(z)/\Psi_\Z|^{p/\lambda} e^{k_\Z(z)}(1 - |z|^2) \le C$. Then $g_0 =
(g/\Psi_\Z)$ will satisfy the appropriate growth with
$\eps=1-\lambda$.

Now for the details.

\begin{lemma}\label{lem:main1}
  Let $0 < p  < \infty$ and suppose that $\Z$ is an interpolating
  sequence for $A^p$. There exist positive constants $C$, $\delta$,
  $\eta$ and $\eps$ depending only on $p$ and the interpolation
  constant $M$ of $\Z$, and a family of analytic functions $\{ g_a \st a
  \in \ID\}$ such that $| g_a e^{k_\Z}| > \delta$ in $D(a,\eta)$ and
  $|g_a(z)e^{k_\Z(z)}|^p(1 - |M_a(z)|^2)^{1-\eps} \le C$ for all
  $z\in \ID$.
\end{lemma}

\begin{proof}
Let us first assume that $D(0,r)$ is disjoint from
$\Z$. If we choose $\lambda \in (0,1)$ close enough to $1$ then $\Z =
p_\lambda(\Z')$ for some sequence $\Z'$ in $\ID$. For each $a \in \Z$
let $a' \in \Z'$ satisfy $\arg a' = \arg a$ and $1-|a|^2 = \lambda(1 -
|a'|^2)$, that is $a = p_\lambda(a')$. It is straightforward to see that
$\sup_{a\in \Z} \psi(a,a')$ tends to $0$ as $\lambda$ tends to $1$.
Therefore, by Proposition~\ref{prop:interpolationperturb}, $\Z'$ will be
an interpolation sequence for $A^p$ if $\lambda$ is sufficiently close
to $1$ ($\lambda$ depending only on $r$, $p$ and $M$), with an
interpolation constant close to $M$.

By Corollary~\ref{cor:lowerestimate} there is a $\delta_1 > 0$ depending
only on $r$, $M$ and $p$, and a function $f$ in the unit ball of $A^p$
which has $\Z(f) = \Z'$ and $|f(0)| > \delta_1$. Then by
Proposition~\ref{prop:zeroperturb} we can find a function $h$ in the
unit ball of $A^{p/\lambda}$ with $\Z(h) = \Z$ and $|h(0)| > \delta_2 =
\delta_1^\lambda/C_1$ with $C_1$ depending only on $\lambda$ and the
separation constant of $\Z'$ (which ultimately depends only on $p$ and
$M$).

By Proposition~\ref{prop:growth} with $p$ replaced by $p/\lambda$, we
obtain a zero-free function $f^*$ satisfying $| f^*(0)| > \delta_2$ and
the growth condition $|f^*(0)|^{p/\lambda}(1 - |z|^2) < C_2$. This
growth condition imposes a bound on the absolute value of the derivative
of $f^*$ on compact sets, and so we have $|f^*(z)| > \delta_3 =
\delta_2/2$ in some neighborhood $D(0,\eta)$ of the origin. By the
second part of Proposition~\ref{prop:growth} we obtain $g_0$ satisfying
\begin{gather}
    \left| g_0(z) e^{k_\Z(z)} \right|^{p/\lambda} (1 - |z|^2) < C_3,\quad
    z\in \ID \label{eq:necbd1}\\
\intertext{and}
    |g_0(z) e^{k_\Z(z)}| > \delta_2, \quad z \in D(0,\eta).
    \label{eq:necbd2}
\end{gather}
The second inequality follows because $|g_0| e^{k_\Z} = |f^*|/\sigma$,
and $\sigma(z) < 1$ for all $z\in \ID$.

Now let $\Z$ be arbitrary and let $r$ be the separation constant of
$\Z$. Then $D(0,r)$ contains at most one point of $\Z$. If there is such
a point and we let it be $b$, apply the above argument to $\Z\setminus
\{ b \}$ to obtain $g_0$. This $g_0$ satisfies the stated conditions for
$\Z \setminus \{ b \}$, but it also satisfies them for the full sequence
$\Z$ with $C_3$ and $\delta_2$ replaced by slightly different $C_4$ and
$\delta$. This is because because the weight functions for the different
sequences differ only by the single factor $\exp \left( \frac {|z|^2}
{2} \frac {(1 - |b|^2)^2} {| 1 - \bar b z |^2} \right)$, which is
bounded above and away from 0, the bounds depending only on $r$.

To get $g_a$ for arbitrary $a$ consider the interpolating sequence $\Z_a
= M_a(\Z)$ obtained by subjecting $\Z$ to the M\"obius transformation
$M_a$. This has the same interpolation constant as $\Z$. Moreover, the
function $k_{\Z_a}$ differs from $k_{\Z}\circ M_a$ by a harmonic
function. This follows from the easily calculated formula
\begin{equation}\label{eq:Moebiusinv}
    (1 - |z|^2)^2\lap k_Z = \sum_{b\in \Z} \left( 1 - \left| M_b(z)
    \right|^2 \right)^2 = \sum_{b\in \Z} (1 - \psi(z,b)^2)^2,
\end{equation}
and the M\"obius invariance of both sides,
Therefore, if we obtain a $g_0$ as above for the sequence $\Z_a$ we have
\begin{equation}
    |g_0e^{k_{\Z_a}}| = |g_0|e^u e^{k_\Z \circ M_a}.
\end{equation}
for some harmonic function $u$. If we now let $h = u + iv$ be analytic,
and let $g_a = g_0 (M_a) \exp h(M_a)$, then from
\eqref{eq:necbd1} and \eqref{eq:necbd2} for $\Z_a$ we get for $g_a$ the
estimates
\begin{gather}
    g_a(z)e^{k_\Z(z)} > \delta > 0, \quad z\in D(a,\eta)\\
\intertext{and}
    |g_a e^{k_\Z}|^{p}(1 - |M_a(z)|^2)^\lambda \le C_4^\lambda
\end{gather}
Which is what was required, with $\eps = 1 - \lambda$ and
$C=C_4^\lambda$.
\qed\end{proof}

Now we are ready to construct a solution kernel for the $\dbar$-equation
on the weighted space $L_\Z^p$.

\begin{theorem}\label{thm:main}
If $\Z$ is an interpolating sequence for $A^p$, $p \ge 1$, then there
exists covering of $\ID$ by disks $D(a_j,\eta)$, a partition of unity
$\beta_j$ subordinate to this covering, functions $g_{a_j}$ as in
Lemma~\ref{lem:main1}, and a positive number $m$ such that
\begin{equation}
  u(z) = \frac{1}{\pi} \sum_{j=1}^\infty g_{a_j}(z) \int_{\ID}
    \frac {\beta_j(w) f(w)} {g_{a_j}(w)} \frac {(1 - |w|^2)^{m-1}}
    {(z - w) ( 1 - \bar wz )^m} \,dA(w)
\end{equation}
satisfies the equation $(1 - |z|^2)\dbar u = f$ and the estimate $\| u
\|_{p,\Z} \le C \| f \|_{p,\Z}$ for some constant $C$ independent of $f$.
\end{theorem}

\begin{proof}
We obtain in the usual way a covering of $\ID$ by disks
$D_j=D(a_j,\eta)$ which also satisfiies a finite overlap condition
$\sup_{z\in \ID} \sum \chi_{D_j}(z) < \infty$ by taking a maximal
sequence $a_j$ such that $D(a_j, \eta/2)$ is disjoint.

We observe that if $m > 0$ then each individual term $u_j$ in the
sum defining $u$ satisfies $(1 - |z|^2)\dbar u_j = \beta_j f$. It is
enough to know that
\begin{equation*}
    v(z) = \frac {1} {\pi} \int \phi(w) \frac{(1 - |w|^2)^m}
        {(z - w) (1 - \bar w z)^m} \,dA(w)
\end{equation*}
satisfies $\dbar v = \phi$ for any compactly supported $\phi$. This
follows from the well-known property of the kernel $(z - w)^{-1}$ and the
fact that this kernel differs from that one by a function analytic in
$z$.

Now we want to prove the estimates. Multiplying both sides by
$e^{k_\Z(z)}$, and rewriting the integrals in terms of $f_1(w) =
f(w)e^{k_\Z(w)} \in L^p$, it is enough to show that the kernel
\begin{equation}\label{eq:solutionkernel}
  \frac{1}{\pi} \sum_{j=1}^\infty \left| g_{a_j}(z) e^{k_\Z(z)} \frac
    {\beta_j(w)} {g_{a_j}(w)e^{k_\Z(w)}} \frac {(1 - |w|^2)^{m-1}} {(z -
    w) ( 1 - \bar wz )^m}  \right|
\end{equation}
defines a bounded integral operator from $L^p(dA(w))$ to itself. The
lower bounds in Lemma~\ref{lem:main1} give $g_{a_j}(w)e^{k_\Z(w)} >
\delta$ in the support of $\beta_j$, and the upper bounds imply that
\begin{equation}
    |g_{a_j}(z)e^{k_\Z(z)}|^p \le C(1 - |M_{a_j}(z)|^2)^{\eps - 1} = C
    \frac { | 1 - \bar a_j z |^{2-2\eps} } {(1 - |a_j|^2)^{1-\eps}
    (1 - |z|^2)^{1-\eps} }
\end{equation}
for some positive $\eps$ and $C$. Moreover, for $w$ in the support
of $\beta_j$ we have $1 - |a_j|^2 \ge c(1-|w|^2)$ and $|1 - \bar a_j z|
\le C |1 - \bar w z|$ for all $z \in \ID$. Using these estimates in
\eqref{eq:solutionkernel}, we need only prove estimates for the kernel
\begin{multline}
    \sum_j  \frac { | 1 - \bar a_j z |^{(2-2\eps)/p} \beta_j(w) } {
        ( 1 - |a_j|^2 )^{(1-\eps)/p} ( 1 - |z|^2 )^{(1-\eps)/p} }
        \frac { ( 1 - |w|^2 )^{m-1} } { | 1 - \bar wz |^m }  \\
    {}\le C \sum_j \beta_j(w) \frac{ ( 1 - |z|^2 )^{(\eps-1)/p}
        (1 - |w|^2)^{m-1+(\eps-1)/p}} { |z-w| \, | 1 - \bar wz
        |^{m+2(\eps-1)/p} } \\
    =   \frac { ( 1 - |z|^2 )^{(\eps-1)/p} ( 1 - |w|^2
        )^{m-1+(\eps-1)/p} } { |z-w| \, | 1 - \bar wz
        |^{m+2(\eps-1)/p} }
\end{multline}
Now we invoke Lemma~\ref{lem:boundedkernel} with $a =
(\eps-1)/p $ and $b = m - 1 + (\eps -1)/p$. Clearly $a > -1/p$
and we will also have $b > -1/p'$ provided we choose $m \ge 2/p$. In
particular $m = 2$ suffices for all $p \ge 1$.

Finally, we need to show that $u$ satisfies the $\dbar$-equation. Fix $z
\in \ID$ and note that only finitely many $\beta_j$ have support
intersecting $D(z,\eta)$. Based on the above estimates, it is
straightforward to see that in this disk, the sum over the remaining
terms is a convergent sum of functions analytic in $D(z,\eta)$ and so $u$
differs from an analytic function by a finite sum $\sum u_j$ with  each
term satisfying $(1 - |z|^2)\dbar u_j(z) = \beta_j(z)f(z)$.
Thus, $(1 - |z|^2)\dbar u(z) = \sum \beta_j(z) f(z) = f(z)$, as
required.
\qed\end{proof}

\begin{remark}
While the functions $g_a$ used in constructing this solution operator
were zero-free, there is really no need for this in the proof. We needed
the precursors of $g_a$ in Proposition~\ref{prop:zeroperturb} to be
zero-free so we could take fractional powers, but all we need here was
the local lower bound and a global upper bound estimate. We return to
this observation in section~\ref{sec:ortega}.
\end{remark}

\section{When $p < 1$}\label{sec:plt1}

When $p < 1$, the proof of Theorem~\ref{thm:main} cannot hold as it
stands because the integrals defining $u$ need not exist if $f$ is only
assumed to be in $L^p$. However, it is possible to replace $L^p_\Z$ with
a space that is locally like $L^1_\Z$ but globally like $L^p_\Z$ and
which contains $A^p_\Z$.

To this end we fix a number $R$ in the interval $(0,1)$ and let $D_z$
denote the pseudohyperbolic disk with center $z$ and pseudohyperbolic
radius $R$. For a function $f$ with $|f|^q$ locally integrable, we
define the \term{local means}, $m_q(f)$ by
\begin{equation*}
    m_q(f,z) = \left(\frac {1} {|D_z|} \int_{D_z} |f(w)|^q \,dA(w)
        \right)^{1/q}.
\end{equation*}
Then let $L^p_q$ denote the set of measurable functions $f$ such that
$m_q(f)$ belongs to $L^p$ and let $L^p_{q,\Z}$ denote the functions $f$
such that $fe^{k_\Z}$ belongs to $L^p_q$. Let $A^p_q$ and $A^p_{q,\Z}$
denote the subspaces consisting of all analytic functions in the
respective space. The space $A^p_q$ is actually the same as $A^p$.
Moreover, as long as $\Z$ is a separated sequence, $A^p_{q,\Z}$ is the
same as $A^p_\Z$, though we will not need this fact. We will use $\| f
\|_{p,q}$ to denote the $L^p(dA)$ norm of $m_q(f)$ and $\| f
\|_{p,q,\Z}$ will denote the $L^p$ norm of $m_q(fe^{k_\Z})$.

\begin{theorem}\label{thm:mainplt1}
If $p < 1$ and $\Z$ is separated in the pseudo-hyperbolic metric, then
the following are equivalent:
\begin{enumerate}
\renewcommand{\theenumi}{\alph{enumi}}
\renewcommand{\labelenumi}{\textup{(\theenumi)}}
\item \label{first-item}$\Z$ is an interpolating sequence for $A^p$.
\item For every $q \ge 1$ there is a constant $C$ such that the
    $\dbar$-equation \eqref{eq:dbar} has for any $f \in
    L^p_{q,\Z}$ a solution $u$ satisfying $\| u \|_{p,\Z} \le C\| f
    \|_{p,q,\Z}$.
\item \label{third-item}There exists a $q \ge 1$ and a constant $C$ such
    that the $\dbar$-equation \eqref{eq:dbar} has for any $f
    \in L^p_{q,\Z}$ a solution $u$ satisfying $\| u \|_{p,\Z} \le C\| f
    \|_{p,q,\Z}$.
\end{enumerate}
\end{theorem}

\begin{proof}
Suppose (\ref{third-item}) holds, then the function $(1 -
|z|^2)g(z)/\Psi(z)$ constructed in Theorem~\ref{thm:easydirection} also
belongs to $L^p_{q,\Z}$ so the identical proof and the fact that $A^p_q
= A^p$ gives us (\ref{first-item}).

Now suppose $\Z$ is an interpolating sequence for $A^p$ with $p < 1$, and
Let $q\ge 1$ be arbitrary and let $f \in L^p_{q,\Z}$.

The existence of the family $g_a$ used in Theorem~\ref{thm:main} was
shown for every $p > 0$. Moreover, the integrals making up the solution
operator constructed there are defined because $fe^{k_\Z}$ is locally in
$L^q$ with $q \ge 1$. It therefore suffices to show that it satisfies
the appropriate bounds. To this end we only need the appropriate
generalization of Lemma~\ref{lem:boundedkernel}, which is
Lemma~\ref{lem:schurplt1} below.
\qed\end{proof}

When estimating bounds on the types of integral operators we encounter
here, a sort of meta-theorem is that best results are obtained for $p >
1$ using Schur methods, that is, using  Lemma~\ref{lem:schur} or the
techniques of its proof; while for $p < 1$, best results seem to come
from discretizing the integrals involved and applying the inequality
between the $p$th power of a sum and the sum of $p$th powers. In this
case, where $p < 1$ but $q \ge 1$, we will need to use both.

First we need a lemma on discretizing the $L^p_q$ norm.

\begin{lemma}\label{lem:discrete}
Let $\{ z_k \}$ be a sequence in $\ID$ with reparation constant at least
$R/2$ and such that the disks $D(z_k,R/2)$ cover $\ID$. Then for any
$0 < p < \infty$ there are constants $c_1$ and $C_2$ such that
\begin{equation}\label{eq:contdiscrete}
   c_1 \sum_k (1 - |z_k|^2)^2 m_q(f,z_k)^p \le \int m_q(f)^p \,dA \le
   C_2 \sum_k (1 - |z_k|^2)^2 m_q(f,z_k)^p
\end{equation}
Moreover, the space $L^p_q$ is independent of the radius $R$ used in its
definition.
\end{lemma}

\begin{proof}
Let us temporarily show the dependence of $m_q$ on $R$ (and suppress the
dependencs on $q$) by writing $m_R$ in place of $m_q$. Let $R_- < R < R^+$ be
defined as follows: $R^+$ is the radius of the disk formed by taking the
union of $D(z,R)$ over all $z \in D(0,R)$, In formula: $R^+ = 2R/(1 +
R^2)$. And $R_-$ is chosen so that $R = 2R_-/(1 + R_-^2)$. It follows
that $R_- > R/2$.

Since the disks $D(z_k, R/4)$ are disjoint while the disks $D(z_k, R/2)$
cover $\ID$, we can write the integral above as a sum of integrals over
disjoint sets $E_k$ satisfying $D(z_k, R/4) \subset E_k \subset D(z_k,
R/2)$. Note that if $z \in E_k$ then $E_k \subset D(z, R) \subset D(z_k,
R^+)$. Therefore, on $E_k$ we have $m_R(f,z)^p \le C_R m_{R^+}(f,z_k)^p$.
Integrating this inequality over $E_k$ and summing on $k$ gives the
second inequality in \eqref{eq:contdiscrete} with $m_{R^+}$ used for
$m_q$ on the right side.

We also have, if $z \in E_k$, then $D(z_k, R_-)$ is contained in $D(z,
R)$, so $m_{R_-}(f,z_k)^p \le C_R m_R(f,z)^p$ for all $z \in E_k$.
Integrating this over $E_k$ and summing gives the other inequality, but
for $m_{R_-}$ in place of $m_q$ on the left side.

Finally, each disk $D(z_k, R^+)$ is covered by the finite number
(independent of k) of $D(z_j, R_-)$ that intersect it, and so $\sum (1 -
|z_k|^2)^2 m_{R^+}(f,z_k)^p \le C \sum (1 - |z_k|^2)^2
m_{R_-}(f,z_k)^p$, which shows that the flanking sums are equivalent.

Repeating the process indefinitely with a new $R = R^+$ or $R=R_-$
shows that the discrete sums and the integrals are equivalent norms for
all $R$.
\qed\end{proof}

\begin{lemma}\label{lem:schurplt1}
If $p < 1$ and $q \ge 1$, then the kernels
\begin{gather*}
  K(z,w) =  \frac{(1 - |z|^2)^a(1 - |w|^2)^b}{|1 - \bar w z|^{a+b+2}}\\
  B(z,w) =  \frac{(1 - |z|^2)^a(1 - |w|^2)^b}{|z-w| \, |1 - \bar w
    z|^{a+b+1}}
\end{gather*}
define bounded operators on $L^p_q(dA)$ \textup(integration with respect to
$w$\textup) provided $a > -1/p$, $b > 2/p - 1/q - 1$ and $a  + b + 2 > q/p $.
\end{lemma}

\begin{proof}
Since $B(z,w) \ge K(z,w)$, it suffices to prove this for the kernel $B$
only. Let $R$ be some convenient radius, and select a sequence $\{z_k\}$
as in Lemma~\ref{lem:discrete}. Let $D_k$ denote $D(z_k,R)$ and recall
that  $1 - |w|^2$ is equivalent to $1 - |z_k|^2$ when $w \in D_k$, and
$|1 -\bar w z|$ is equivalent to $|1 -\bar z_k z|$, and $|D_k|$ is
equivalent to $(1 - |z_k|^2)^2$. Therefore,
\begin{align*}
    Bf(z)
        &\defeq \int B(z,w) f(w) \,dA(w)\\
        &\le  \sum_{k} \int_{D_k} K(z,w) |f(w)| \left| \frac{1 - \bar wz}
            {z-w} \right| \,dA(w) \\
        &\le  \sum_k K(z,z_k)(1 - |z_k|^2)^2 \frac{1}{|D_k|}
            \int_{D_k} |f(w)| \left| \frac{1 - \bar w z}{z - w} \right|
            \,dA(w)\\
        &=  \sum_{k} K(z,z_k)(1 - |z_k|^2)^2 b_k(z)
\end{align*}
where, for notational convenience, we have put
\begin{equation*}
   b_k(z) = \frac{1}{|D_k|} \int_{D_k} |f(w)|
             \left| \frac{1 - \bar w z}{z - w} \right| \,dA(w).
\end{equation*}
Applying H\"older's inequality to the sum, we get
\begin{equation}
   Bf(z)^q
        \le \left(\sum K(z,z_k) (1 - |z_k|^2)^{2 + \alpha q} b_k(z)^q
            \right) \left( \sum_k K(z,z_k)(1 - |z_k|^2)^{2 - \alpha q'}
            \right)^{q/q'}
\end{equation}
The second sum is less than a constant times $(1 - |z|^2)^{-\alpha q'}$
by the discrete version of inequality \eqref{eq:Forelli-Rudin}. (See
\cite{Lue99c} Lemma~3 for a general version.)  This inequality requires
$\alpha$ to be chosen so that $1 < b + 2 - \alpha q' < a + b + 2$.
Putting it in the above and integrating the result over $D_j$ with
respect to $z$ we get
\begin{equation}\label{eq:Bfestimate}
    m_q(Bf,z_j)^q
        \le \sum _k \frac { (1 - |z_j|^2)^{a - \alpha q} (1 -
            |z_k|^2)^{2 + \alpha q} } { | 1 - \bar z_k z_j|^{a+b+2} }
            m_q(b_k,z_j)^q.
\end{equation}

When the distance from $D_j$ to $D_k$ is greater than $R$ the factor
multiplying $|f(w)|$ in the definition of $b_k(z)$ is at most $1/R$, so
$b_k(z) \le C |D_k|^{-2} \int_{D_k} |f| \,dA \le C m_q(f,z_k)$, and
therefore $m_q(b_k,z_j) \le C m_q(f,z_k)$.

On the other hand when the distance from $D_j$ to $D_k$ is less than $R$
and $z\in D_j$, then $|1 - \bar wz|$ is equivalent to $1 - |z_k|^2$ and
so $b_k(z)$ is bounded by a constant times the convolution of
$|f|\chi_{D_k}$ with $(1 - |z_k|^2)^{-1}|z|^{-1}\chi_{S_k}$, where $S_k$ is
a disk centered at $0$ with radius equal to the Euclidean diameter of
$D_j \cup D_k$, which is $C(1 - |z_k|^2)$. Using Young's inequality, we
have
\begin{equation}
  m_q(b_k,z_j) \le \frac{C}{|D_j|^{1/q}} \| \,|f|\chi_{D_k} \|_q \|(1 -
    |z_k|^2)^{-1} |z|^{-1}\chi_{S_k} \|_1 \le C m_q(f,z_k).
\end{equation}
Putting these estimates for $m_q(b_k,z_j)$ into \eqref{eq:Bfestimate} we
arrive at
\begin{equation}
    m_q(Bf,z_j)^q
        \le \sum _k \frac { (1 - |z_j|^2)^{a - \alpha q} (1 -
            |z_k|^2)^{2 + \alpha q} } { | 1 - \bar z_k z_j|^{a+b+2} }
            m_q(f,z_k)^q.
\end{equation}
Now we need to raise this to the $p/q$ power and bring the power inside
the summation
\begin{equation}
    m_q(Bf,z_j)^p
        \le C \sum _k \frac { (1 - |z_j|^2)^{(a - \alpha q)p/q} (1 -
            |z_k|^2)^{(2 + \alpha q)p/q } } { | 1 - \bar z_k z_j
            |^{(a+b+2)p/q} } m_q(f,z_k)^p,
\end{equation}
and, finally, multiply by $(1 - |z_j|^2)^2$ and sum on $j$.
\begin{align*}
    \| Bf \|_{p,q}^p
        &\le C\sum (1 - |z_j|^2)^2 m_q(Bf,z_j)^p\\
        &\le C \sum _k \sum_j \frac { (1 - |z_j|^2)^{(a - \alpha q)p/q +
            2} (1 - |z_k|^2)^{(2 + \alpha q)p/q } } { | 1 - \bar z_k z_j
            |^{(a+b+2)p/q} } m_q(f,z_k)^p\\
        &\le C\sum_k (1 - |z_k|^2)^2 m_q(f,z_k)^p\\
        &\le C\| f \|_{p,q}^p
\end{align*}
Where the next-to-last inequality again uses the discrete version of
\eqref{eq:Forelli-Rudin} and requires $\alpha $ to be chosen satisfying
$1 < (a - \alpha q)p/q + 2 < (a + b + 2) p/q$.

Thus, everything depends on the existence of a real number $\alpha$
satisfying
\begin{gather*}
    1 < b + 2 - \alpha q' < a + b + 2,\\
\intertext{and}
    1 < (a - \alpha q)p/q + 2 < (a + b + 2) p/q,
\end{gather*}
or equivalently
\begin{gather*}
    -\frac{a}{q'} < \alpha < \frac{b + 1}{q'} \\
\intertext{and}
    2/p - (b + 2)/q < \alpha < a/q + /p.
\end{gather*}
Such an $\alpha$ will exist if each left side is less than both right
sides, leading to the following four requirements:
\begin{gather*}
0 < a + b + 1\\
0 < a + \frac{1}{p}\\
\frac {q}{p} < a + b + 2\\
\frac {2}{p}< b  + \frac {1}{q} + 1 \\
\end{gather*}
The last three of these are the hypotheses of this lemma, and the first
one is implied by the third, so the lemma is proved.
\qed\end{proof}

Returning to Theorem~\ref{thm:mainplt1}, the arguments of
Theorem~\ref{thm:main} lead to the requirement that the kernel
\begin{equation*}
     \frac { ( 1 - |z|^2 )^{(\eps-1)/p} ( 1 - |w|^2
        )^{m-1+(\eps-1)/p} } { |z-w| \, | 1 - \bar wz
        |^{m+2(\eps-1)/p} }
\end{equation*}
define a bounded operator. In the present context we require that
$(\eps - 1)/p > -1/p$, which is clear, and also that $m + 1 +
(\eps - 1)/p > 2/p + 1/q'$ and $m + 1 + 2(\eps - 1)/p > q/p$,
both of which can be made to hold by choosing $m$ sufficiently large.

\section{Seip's density condition}\label{sec:seip}

Here we show that Seip's density condition (\cite{Sei94}) for
interpolation is equivalent to a condition on $k_\Z$. First we recall
the definitions and the condition.

Let $\Z$ be a sequence in $\ID$ separated in the hyperbolic metric. For
$b\in \ID$ let $\Z_b = M_b(\Z)$. We define $D^+(\Z,r)$ by
\begin{equation*}
    D^+(\Z,r) =  \frac{\displaystyle \sup_{b\in \ID} \sum_{a\in\Z_b,
        1/2<|a|<r} \log \frac {1} {|a|}} {\displaystyle\log\frac {1} {1
        - r}}.
\end{equation*}
then the \term{upper uniform density} of $\Z$ is $D^+(\Z) = \limsup_{r
\to 1} D^+(\Z, r)$. Given that the denominator is unbounded, we may
replace the numerator by any quantity that differs from the given one by
a constant independent of $b$ and $r$, we choose to use
\begin{equation*}
    D^+(\Z) = \limsup_{r\to 1}  \frac {\displaystyle \sup_{b\in \ID}
       \sum_{a\in \Z_b, |a|<r} \frac {1 - |a|^2} {2}}
       {\displaystyle \log \frac {1} {1 - r^2}}.
\end{equation*}
Then Seip's density condition for interpolation in $A^p$ is $D^+(\Z) <
1/p$.

We now show that the numerator can be replaced by
\begin{equation}
   \frac {1} {2\pi} \int_0^{2\pi} k_\Z (r e^{i\theta}) \,d\theta = \frac
        {r^2} {2} \sum_{a\in\Z} \frac {(1 - |a|^2)^2} {1 - |a|^2r^2}.
\end{equation}
Let us define the alternative density
\begin{equation}
    S^+(\Z) = \limsup_{r\to 1}  \frac { \displaystyle \sup_{b\in \ID}
    \int_0^{2\pi} k_{M_b(\Z)} (r e^{i\theta}) \,d\theta} {\displaystyle
    \log \frac {1} {1 - r^2}}.
\end{equation}
Obviously this remains unchanged if we replace the $r$ with $r_\eps$ where
$1-r_\eps^2 = \eps(1 - r^2)$. But since $\log(1 - r_\eps^2) =
\eps + \log(1 - r^2)$, we get the same $\limsup$  using $r_\eps$
in the numerator and just $r$ in the denominator. Let us now estimate
the difference in the numerators.
\begin{equation*}
   \sum_{|a| < r} \frac{1 - |a|^2}{2} - \frac {1} {2\pi} \int_0^{2\pi}
   k_\Z (r_\eps e^{i\theta}) \,d\theta.
\end{equation*}
Clearly it equals
\begin{equation}\label{eq:comparison}
    \sum_{|a| < r} \frac{1 - |a|^2}{2} \left( 1 - \frac{(1 -
    |a|^2)r_\eps^2}{1 - |a|^2r_\eps^2} \right) - \frac {r_\eps^2} {2}
    \sum_{|a| \ge r} \frac{(1 - |a|^2)^2}{1 - |a|^2r_\eps^2}
\end{equation}
The first sum simplifies to
\begin{equation*}
    \sum_{|a| < r} \frac{1 - |a|^2}{2} \left( \frac {(1 - r_\eps^2)}
    {1 - |a|^2r_\eps^2} \right)
\end{equation*}
and since
\begin{equation*}
    \frac {(1 - r_\eps^2)} {1 - |a|^2 r_\eps^2} \le \frac {(1 -
    r_\eps^2)} {1 - |a|^2} < \frac {(1 - r_\eps^2)} {1 - r^2} = \eps
\end{equation*}
when $|a| < r$, we have
\begin{equation}\label{eq:head}
    (1 - \eps) \sum_{|a| < r} \frac {1 - |a|^2} {2} \le \frac {r_\eps^2}
        {2} \sum_{|a| < r} \frac {(1 - |a|^2)^2} {1 - |a|^2 r_\eps^2}
        \le \sum_{|a| < r} \frac {1 - |a|^2} {2}
\end{equation}
As for the second sum in \eqref{eq:comparison}, we first observe that
$\sum_{|a| \ge r} (1 - |a|^2)^2$ is dominated by the area of the annulus
$r < |z| < 1$ up to a constant that depends only on the separation
constant of the sequence so
\begin{equation} \label{eq:tail}
    \frac {r_\eps^2} {2} \sum_{|a|\ge r} \frac {(1 - |a|^2)^2} {1 - |a|^2
        r_\eps^2} \le \sum_{|a|\ge r} \frac {(1 - |a|^2)^2} {1 - r_\eps^2}
        \le C \frac {1 - r^2} {1 - r_\eps^2} \le \frac {C} {\eps}.
\end{equation}
These estimates continue to hold for if we replace $\Z$ with $M_b(\Z)$
and so they hold if we take the supremum over $b$.
Combining \eqref{eq:head} and \eqref{eq:tail} with \eqref{eq:comparison}
(for all $M_b(\Z)$), we see that we have $(1 - \eps)D^+(\Z) \le S^+(\Z)
\le D^+(\Z)$. Since $\eps$ is arbitrary, the two densities are the same.

Let us restate the criterion in terms of $S^+$ but in a form not using
$\sup$ or $\limsup$.

\begin{proposition}
In order that $\Z$ be an interpolating sequence for $A^p$ it is
necessary and sufficient that there exist an $\eps > 0$ and $r^* < 1$ such
that for all $w\in \ID$ and all $r$ in the interval $(r_*, 1)$ we have
\begin{equation}\label{eq:means}
    \frac {p} {2\pi} \int_0^{2\pi} k_{M_w(\Z)} (r e^{i\theta}) \,d\theta
        < (1 - \eps) \log \frac {1} {1 - r^2}
\end{equation}
\end{proposition}

We now show that this condition can be replaced with one involving the
Laplacian of $k_\Z$ (easy) and also that it is sufficient to have the
inequality for a single value of $r$ (not as easy).

\begin{corollary}
In order that $\Z$ be an interpolating sequence for $A^p$ it is
necessary and sufficient that it be separted in the pseudo-hyperbolic
metric and there exist an $\eps > 0$ and an $r_*$ with $0 <
r_* < 1$ such
that for all $w\in \ID$ have
\begin{equation}\label{eq:laplacemeans}
    p \int_{|z|<r_*} \lap k_{M_w(\Z)} (z) \log \frac {r_*^2} {|z|^2}
        \,dA(z) < (1 - \eps) \int_{|z|< r_*} \frac {1} {(1 - |z|^2)^2}
        \log \frac {r_*^2} {|z|^2} \,dA(z).
\end{equation}
\end{corollary}

\begin{proof}
The equivalence of inequality \eqref{eq:means} (for a fixed $r=r_*$) and
\eqref{eq:laplacemeans} is just a matter of Green's formula. What we
need to do is show that the \eqref{eq:laplacemeans} for a single value
of $r_*$ implies it for all $r$ sufficiently close to $1$.

We replace $k_{M_w(\Z)}$ with $k_\Z \circ M_w$, which has the same
Laplacian, and then change variables, rewriting \eqref{eq:laplacemeans}
as
\begin{multline*}
    p \int_{D(w,r_*)} \lap k_{\Z}(z) \log \frac
        {r_*^2} {|M_w(z)|^2} \,dA(z)\\ < (1 - \eps) \int_{D(w,r_*)} \frac {1}
        {(1 - |z|^2)^2}\log \frac {r_*^2} {|M_w(z)|^2} \,dA(z).
\end{multline*}
Now we let $r_* < R < 1$ and then integrate the above inequality with
respect to $\log(R^2/|w|^2) \,d\lambda(w)$ over $D(0,R)$, where
$d\lambda(w) = dA(w)/(1 - |w|^2)^2$ is the M\"obius invariant measure on
$\ID$. After exchanging order of integration, we get
\begin{multline}\label{eq:multiline}
    \frac {p} {\pi} \int_{\ID} \int_{\ID} \chi_{D(z,r_*)\cap D(0,R)}(z)
         \log\frac{R^2}{|w|^2}\log \frac {r_*^2} {|M_z(w)|^2}
         \,d\lambda(w) \lap k_{\Z}(z) \,dA(z)\\
        <   (1 - \eps) \int_{\ID} \int_{\ID} \chi_{D(z,r_*)\cap
         D(0,R)}(z) \log \frac {R^2} {|w|^2} \log \frac {r_*^2} {|M_z(w)|^2}
         \,d\lambda(w) \frac{dA(z)}{(1 - |z|^2)^2}
\end{multline}
We remark for later use, that $\lap k_\Z$ is bounded by a constant
(depending only on the separation constant) times $(1 - |z|^2)^{-2}$,
which follows from the formula \eqref{eq:Moebiusinv}.

We split the outermost integrals above into three parts:
\begin{itemize}
\item[(I)] The integral over those $z$ with $|z| \ge r_*$ and $D(z,r_*)
        \subset D(0,R)$. Then $D(z,r_*) \cap D(0,R) = D(z,r_*)$ and
        $\log (R^2/|w|^2)$ is harmonic on $D(z, r_*)$.
\item[(II)] The integral over $|z| < r_*$.
\item[(III)] The integral over those $z$ such that $D(z,r_*)$ meets
        $D(0,R)$ but is not contained in it.
\end{itemize}
In the integral (I), the inner integral on both sides of inequality
\eqref{eq:multiline} is
\begin{equation}\label{eq:inner}
    \int_{D(z,r_*)} \log \frac {R^2} {|w|^2} d\lambda(w) =
        \int_{D(0,r_*)} \log \frac {R^2} {|M_z(w)|^2} \,d\lambda(w) = C
        \log \frac {R^2} {|z|^2}
\end{equation}
where $C$ depends only on $r_*$.

The integral (II) (on both sides of the inequality) is easily estimated
to be bounded by a constant depending only on $r_*$ and the separation
constant of $\Z$.

Finally, for the integral (III), the inner integral is dominated by by
$C(1-R^2)$. And the function $\lap k_\Z$ is less than $C/(1 - |z|^2)^2$
and the set in question is $R_1 < |z| < R_2$ with $1 -R_j^2 \sim (1 -
R^2)$, so that integral (for both the left and right side) is bounded by
a constant depending only on $r_*$ and the separation constant of $\Z$.

Putting these estimates into \eqref{eq:multiline},
we obtain a constant, depending only on $r_*$ and the separation
constant of $\Z$ such that
\begin{equation*}
    p \int_{|z| < R} \lap k_\Z(z) \log \frac
        {R^2} {|z|^2} \,dA(z) < C + (1 - \eps) \int_{|z|< R} \frac {1}
        {(1 - |z|^2)^2}\log \frac {R^2} {|z|^2} \,dA(z).
\end{equation*}
Since the integral on the right grows like $\log(1/(1 - R^2))$ (in fact
it is a constant multiple of it), we can reduce $\eps$ slightly and omit
the $C$ for sufficiently large $R$. Finally, we get
\begin{equation*}
    p \int_{|z| < R} \lap k_{M_w(\Z)}(z) \log \frac
        {R^2} {|z|^2} \,dA(z) < (1 - \eps) \int_{|z|< R} \frac {1}
        {(1 - |z|^2)^2}\log \frac {R^2} {|z|^2} \,dA(z).
\end{equation*}
for all $M_b$ and all sufficiently large $R$ because the estimates
depended only on the $r_*$ and the separation constant of $\Z$.
\qed\end{proof}

For reference below, we can rewrite \eqref{eq:laplacemeans} as
\begin{equation}\label{eq:finaldensity}
    \int_{D(w)} \left(1 - p(1 - |z|^2)^2 \lap
        k_\Z(z) \right) \log \frac {R^2} {|M_w(z)|^2} \,d\lambda(z) >
        \eps \int_{D(w)} \log \frac {R^2} {|M_w(z)|^2} \,d\lambda(z).
\end{equation}
where $D(w) = D(w,r_*)$.

\section{Connection with other weighted $\dbar$-estimates}\label{sec:ortega}

Here we show that a recent result of J.~Ortega-Cerd\`a (\cite{Ort99}) on
weighted $\dbar$-estimates, together with Seip's criterion, implies that
the $\dbar$-equation has solutions with bounds in $L^p_\Z$.
Ortega-Cerd\`a's result deals with the spaces $L^p(\omega \,dA)$ where
$\omega(z) = e^{-\phi(z)}/(1-|z|^2)$. If we want to see the spaces
$L^p_\Z$ into this context, we must put
\begin{equation}\label{eq:phiforkZ}
    \phi(z) = \log \frac{1}{1-|z|^2} - pk_\Z(z).
\end{equation}

The main result in \cite{Ort99} required that $\phi$ be subharmonic,
which need not be the case for \eqref{eq:phiforkZ} (and cannot be for $p
> 1$). So $\lap \phi$ need not be a positive measure.
Another requirement was
that the measure $\lap \phi$ be \term{locally doubling}.
Finally, the main requirement was that there exists
$\eps>0$ and $0<r_*<1$ such that for all $b\in \ID$
\begin{equation}\label{eq:laplower}
   \int_{D(b,r_*)} \lap \phi(z) \,dA(z) \ge \eps
\end{equation}
The conclusion of his result is that the $\dbar$ equation
\eqref{eq:dbar} always has solutions with bounds in the $L^p(\omega
\,dA)$ norm.

We have seen \eqref{eq:finaldensity} that Seip's condition is equivalent
to the following, where $\phi$ is given by \eqref{eq:phiforkZ}: there
exists $\eps>0$ and $0<r_*<1$ such that for all $b \in \ID$
\begin{equation}\label{eq:lapkZlower}
   \int_{D(b,r_*)} \lap \phi(z) \log \frac {r_*^2} {|M_b(z)|^2} \,dA(z)
      \ge \eps
\end{equation}

The similarity between \eqref{eq:laplower} and \eqref{eq:lapkZlower} is
striking. The $\phi$ in the latter inequality need not be subharmonic,
but if it \emph{were} subharmonic, it can be shown that then these two
conditions would be equivalent, and Seip's criterion plus
Ortega-Cerd\`a's result would then imply Theorem~\ref{thm:main}
in section~\ref{sec:main2}.

Now we show (thanks to J. Ortega-Cerd\'a for pointing me in the right
direction) that $\phi$ can be modified in such a way that the difference
between $\phi$ and its modification is bounded, and such that the new
function satisfies the hypotheses of Ortega-Cerd\`a's theorem.

\begin{proposition}
Given $\phi$ as above and $r_*$ as in \eqref{eq:finaldensity}, define
\begin{equation*}
  \phi_*(w) = \frac{1}{\pi \log[1/(1 - r_*^2)]} \int_{D(w,r_*)} \phi(z) \log
        \frac {r_*^2} {|M_w(z)|^2} \,d\lambda(z)
\end{equation*}
If $\Z$ is an interpolating sequence for $A^p$ then $\phi_*$ is
subharmonic and satisfies  $1 \ge (1 - |z|^2)^2\lap \phi_*(z) >
\epsilon$ for some $\eps > 0$. Moreover $\phi_* - \phi$ is bounded.
\end{proposition}

\begin{proof}
The definition of $\phi_*$ is the \term{invariant convolution}
$\phi*g$ of $\phi$ with the rotationally invariant function
\begin{equation*}
g(z) = \frac{1}{\pi \log[1/(1 - r_*^2)]} \chi_{D(0,r_*)}\log(r_*^2/|z|^2)
\end{equation*}
(see M.~Stoll \cite{Sto94}). While it does not appear to be explicitly
stated in \cite{Sto94}, it follows from the discussion there (page
34--37) that when one of the functions is rotationally invariant,
convolution commutes with the \term{invariant Laplacian}. That is, if
$\invlap f(z) = (1 - |z|^2)^2\lap f(z)$, then $\invlap(\phi*g) =
(\invlap \phi)* g$.

The factor preceding the integral in the definition of $\phi_*$ is a
normalization:
\begin{equation*}
  \pi \log[1/(1 - r_*^2)] = \int_{D(w,r_*)} \phi(z) \log \frac {r_*^2}
    {|M_w(z)|^2} \,d\lambda(z).
\end{equation*}
This follows from the Moebius invariance of the right hand side and a
simple calculation for $w = 0$. Moreover, the left hand side of
\eqref{eq:finaldensity} is, after a normalization, equal to $(\invlap
\phi)*g) = \invlap \phi_*$.  Therefore, \eqref{eq:finaldensity} is
equivalent to $\invlap \phi_* > \epsilon$, which is one of the
conclusions we wanted. The upper bound $\invlap \phi_*(z) \le 1$ is
immediate from the fact that $\invlap \phi = 1 - p\invlap k_\Z < 1$.

The fact that interpolating sequences are separated implies that $(1 -
|z|^2) \grad k_\Z(z)$ is bounded. This is also trivially true of $\log [1
/ (1 - |z|^2)]$, and therefore $(1 - |z|^2) \grad \phi(z)$ is bounded.
This means that there is a constant $C$ such that $|\phi(w) - \phi(z)| <
C$ whenever $z \in D(w,r_*)$. Integrating in the $z$ variable we get
$|\phi(w) - \phi_*(w)| < C$, as required.
\end{proof}

The inequalities we have just seen on $\invlap \phi_*$ trivially imply that
the measure $d\mu(z) = \invlap \phi_*(z) \,d\lambda(z)$ is locally
doubling. Moreover, $\mu(D(z,r)) \ge \epsilon\lambda(D(z,r))$ so all
of Ortega-Cerd\'a's conditions are satisfied. Thus, we get
another proof of the main result for $p \ge 1$. The case $p < 1$ (in a
form such as Theorem~\ref{thm:mainplt1}) is a routine extension.

There is another connection between the results obtained here and those
of Ortega-Cerd\`a. For an arbitrary $\phi$ satisfying the hypotheses in
\cite{Ort99}, a family of analytic functions $g_a$ satisfying
inequalities similar to those required in the proof of
Theorem~\ref{thm:main} (namely $|g_a(z)|^p \omega(z)$ should be bounded
below independent of $a$ in a suitable neighborhood $D(z,\eta)$ of $a$,
and should grow at worst like $(1 - |M_a(z)|^2)^{\eps-1}$) would be
enough to obtain the appropriate weighted $L^p$ estimates. Again, one
need only construct one such function and obtain the whole family using
a M\"obius invariance argument.

In \cite{Ort99} a much more special function is constructed: a
single function $f$ with a specified sequence of zeros, satisfying both
upper and lower bounds on $f(z)e^{-\phi}$. The same techniques used in
\cite{Ort99} to construct $f$ can be somewhat simplified to construct
instead the $g_a$. Unlike the $g_a$ in Theorem~\ref{thm:main}, which are
zero-free, these would necessarily have many zeros. But as we've pointed
out, $g_a$ only needs to be free of zeros near $a$.

The details of the construction of such a family $\{ g_a \}$ are so
similar to those of the construction of $f$ in \cite{Ort99} as not to be
worth repeating here.

\newcommand{\noopsort}[1]{}
\providecommand{\bysame}{\leavevmode\hbox to3em{\hrulefill}\thinspace}
\providecommand{\MR}{\relax\ifhmode\unskip\space\fi MR }
\providecommand{\MRhref}[2]{%
  \href{http://www.ams.org/mathscinet-getitem?mr=#1}{#2}
}
\providecommand{\href}[2]{#2}

\end{document}